\documentclass{amsart}
\usepackage{amsfonts,amssymb,amscd,amsmath,bbm,amsthm,txfonts}
\usepackage{accents}
\usepackage{diagrams} 

\newcommand{\newpageforfinallayout}{\newpage}

\newcounter{tmp} 
\theoremstyle{definition}
\newtheorem{Thm}{Theorem}[section]

\newtheorem{Lem}[Thm]{Lemma}
\newtheorem{Cor}[Thm]{Corollary}
\newtheorem{Def}[Thm]{Definition}
\newtheorem{Asm}[Thm]{Assumption}

\newcommand{\fnto}{\ensuremath{\rightarrow}}  % fuer f: A->B
\DeclareMathOperator{\spl}{split}               % domain
\newcommand{\U}{\ensuremath{\bigcup}}                     % bigunion
\newcommand{\pow}{\ensuremath{\mathfrak{P}}}              % powersetsymbol
\newcommand{\card}[1]{{\ensuremath{\rvert#1\lvert}} }     % cardinality

\newcommand{\al}[1]{\ensuremath{{\aleph_{#1}}} }          % aleph_n
\newcommand{\om}[1]{\ensuremath{{\omega_{#1}}} }          % omega_n

\newcommand{\ho}{\ensuremath{^{\omega}}}                  % hoch omega
\newcommand{\hko}{\ensuremath{^{<\omega}}}                % hoch <omega
\DeclareMathOperator{\rank}{rk}                             % rang
\newcommand{\rz}[1]{\ensuremath{^{#1}}}     % relativ zu
\newcommand{\rzM}{\ensuremath{\rz{M}}}                    % relativ zu M
\DeclareMathOperator{\epd}{epd}

\newcommand{\ZFC}{\ensuremath{\text{ZFC}} }               % ZFC
\newcommand{\ZFCminus}{\ensuremath{\text{ZFC}^-} }        % ZFC-
\newcommand{\ZFCx}{\ensuremath{{\text{ZFC}^\ast}} }        % ZFC*
\newcommand{\ZFCxx}{\ensuremath{{\text{ZFC}^{\ast\ast}}} }        % ZFC**
\newcommand{\fa}{{\ensuremath{\forall}} }
\newcommand{\ex}{{\ensuremath{\exists}} }

\newcommand{\AND}{\ensuremath{\&}}
\newcommand{\OR}{\ensuremath{\vee}}
\newcommand{\THEN}{\ensuremath{\,\rightarrow\,}}

\newcommand{\DEFIFF}{\ensuremath{\,:\leftrightarrow\,} }
\newcommand{\DEFEQ}{\ensuremath{\coloneqq}}
\newcommand{\EQDEF}{\ensuremath{\eqqcolon}}

\newcommand{\bB}{\ensuremath{\mathbbm{B}}} % random algebra
\newcommand{\bC}{\ensuremath{\mathbbm{C}}} % cohen forcing
\newcommand{\esm}{\ensuremath{\prec}}  % elem. submodel
\newcommand{\vD}{\ensuremath{\vDash}}

\newcommand{\qb}{\ensuremath{\text{``}}}
\newcommand{\qe}{\ensuremath{\text{''}}}

\newcommand{\mS}[2]{\ensuremath{{\underaccent{\tilde}{\Sigma}^#1_#2}}}
\newcommand{\mP}[2]{\ensuremath{{\underaccent{\tilde}{\Pi}^#1_#2}}}
\newcommand{\mD}[2]{\ensuremath{{\underaccent{\tilde}{\Delta}^#1_#2}}}

\newcommand{\std}[1]{\ensuremath{\check{#1}}}              % standardname
\newcommand{\forc}{\ensuremath{\Vdash}}

\newcommand{\incomp}{\ensuremath{\perp}}
\newcommand{\comp}{\ensuremath{\parallel}}
\newcommand{\vH}{\ensuremath{[H]}}          % fuer t<G>
\newcommand{\vG}{\ensuremath{[G]}}          % fuer t<G>
\newcommand{\vGP}{\ensuremath{[G_P]}}          % fuer t<G>
\newcommand{\vGQ}{\ensuremath{[G_Q]}}          % fuer t<G>
\newcommand{\n}[1]{\underaccent{\tilde}{#1}}

\newcommand{\hcn}[1]{\n{#1}}

\DeclareMathOperator{\ro}{ro}
\newcommand{\h}{{\hcn{\eta}}}  % name eta des gen reals
\newcommand{\es}{{\eta^\ast}}   % konkretes generic real
\newcommand{\eot}{{\eta^\otimes}}   % konkretes generic real

\newcommand{\BC}{\text{BC}}  % I^dx+
\newcommand{\IB}{I_\text{BC}}  % I^dx+
\newcommand{\IpB}{I^+_\text{BC}}  % I^dx+
\newcommand{\Ip}{I^+}  % I^dx+

\DeclareMathOperator{\Gen}{Gen}        % Gen(M) = Set of M-generic reals
\newcommand{\aGen}{\Gen^{\text{abs}}}  % Genabs(M) = Set of abs M-generic reals

\newcommand{\bR}{\ensuremath{\mathbbm{R}}}

\begin{document}
\date{December 2002}

\title{Preserving Non-Null with Suslin$^+$  forcings}
\author{Jakob Kellner}
\address{Institut f\"ur Diskrete Mathematik und Geometrie\\
 Technische Universit\"at Wien\\
 1050 Wien, Austria}
\email{kellner@fsmat.at}
\urladdr{http://www.logic.univie.ac.at/$\sim$kellner}
\thanks{Partially supported by FWF grant P17627-N12.
I thank a referee for pointing out an error and several 
unclarities.}

%========================================================================
%  ABSTRACT
%========================================================================

\begin{abstract}
We %(i.e. I)
introduce the notion of effective Axiom A and use it to show that
some popular tree forcings are Suslin$^+$.  We introduce transitive nep and 
present a simplified version of Shelah's ``preserving a little implies
preserving much'': If $I$ is a Suslin ccc ideal (e.g. Lebesgue-null or meager)
and $P$ is a transitive nep forcing (e.g. $P$ is Suslin$^+$) and $P$ doesn't
make any $I$-positive {\em Borel} set small, then $P$ doesn't make {\em any}
$I$-positive set small. 
\end{abstract}

\maketitle

\section{Introduction}

Properness is a central notion for countable support iterations: If a forcing
$P$ is proper then it is ``well behaved'' in certain respects (most notably $P$
doesn't collapse $\om1$); and properness is preserved under countable support
iterations. Properness can be defined by the requirement that the generic
filter (over $V$) is generic for a countable elementary submodel $N$ 
as well (see \ref{def:proper}).

It turns out that it can be useful to 
require genericity for non-elementary models $M$ as well.\footnote{For 
  this to make sense the forcing notion $P$ has to be definable, otherwise we
  do not know how to find $P$ in $M$, and therefore cannot formulate that 
  $G$ is $P$-generic over $M$.}
The first notion of this kind was Suslin proper \cite{JdSh:292},
with the important special case Suslin ccc. This notion was generalized 
to Suslin$^+$ \cite{tools}. In this paper we recall these definitions,
and introduce an effective version of Axiom A as a tool to show that 
all the usual Axiom A forcings
%\footnote{Axiom A was presented by
%Baumgartner in \cite{Baumgartner}.} 
are in fact Suslin$^+$.

In \cite{Sh:630} Shelah introduced a further generalization: non-elementary
proper (nep) forcing.  Here, the models $M$ considered are not only
non-elementary but also non-transitive.  This allows to deal with long
forcing-iterations (which can never be element of a transitive countable model),
but this also brings some unpleasant technical difficulties.  To avoid some of
these difficulties \cite{Sh:630} uses a set theory with ordinals as urelements.

In this paper we define a special case, the ``transitive version'', of nep.  In
this version we consider transitive candidates only, which makes the whole
setting much easier.

As an example of how to apply non-elementary properness we give
a simplified proof of Shelah's ``preserving a little implies preserving much''
\cite[sec. 7]{Sh:630}: If a forcing $P$ is provably nep
and provably doesn't make the set of all old reals Lebesgue null,
then $P$ doesn't make any positive set null.
The proof uses the fact that we can find generic conditions for models of the
form $N[G]$, where $N$ is (a transitive collapse of) an elementary submodel and
$G$ an {\em internal} $N$-generic filter (i.e. $G\in V$).

The proof works in fact not only for the ideal of Lebesgue null sets,
but for all Suslin ccc ideals (e.g. the meager ideal). 
A couple of theorems of this kind lead up to the general
case in \cite{Sh:630}: For the meager case the result is due to Goldstern and
Shelah \cite[Lem XVIII.3.11, p.920]{Sh:f}, the Lebesgue null case in the
special case of $P$=Laver was done by Pawlikowski \cite{PawLav} (building on
\cite{JdSh:308}).
The definition and basic properties of Suslin ccc ideals have been used for a
long time, for example 
in works of Judah, Bartoszy\'{n}ski and Ros{\l}anowski, cited
in \cite{tomekbook}; also related is \cite[\S31]{Si64}.

The result is useful for positivity preservation in limit-steps of
countable support proper
iterations $(P_\alpha)_{\alpha<\delta}$: while it is not clear how one could
argue directly that $P_\delta$ still is Borel positivity preserving, the
equivalent ``preservation of generics'' (see definition \ref{def:presgen}) has
a better chance of being iterable.
In section \cite[XVIII.3.10]{Sh:f} this iterability is claimed for $I$=meager.
For $I$=Lebesgue null the result will appear in \cite{KrSh:828}. 

\subsection*{Annotated contents}
\begin{list}{}{\setlength{\leftmargin}{0.5cm}\addtolength{\leftmargin}{\labelwidth}}
%\begin{description}\setlength{\leftmargin}{0cm}
\item[Section~\ref{sec:suslinplus}, p. \pageref{sec:suslinplus}:]
 We will recall the definition and basic properties of Suslin
proper, Suslin ccc and Suslin$^+$ forcings, and introduce the notions
transitive nep and effective Axiom A. 
We use effective Axiom A to show that Laver, Sacks and similar tree forcings are Suslin$^+$.
\item[Section~\ref{sec:suslinccc}, p. \pageref{sec:suslinccc}:]
We introduce Suslin ccc ideals an their basic properties.  Such ideals are
defined by a Suslin ccc forcing $Q$ with a name for a generic real $\n\eta$ in
the same way as Lebesgue null can be defined from random forcing or meager from
Cohen forcing.
\item[Section~\ref{sec:preservation}, p. \pageref{sec:preservation}:]
We prove Shelah's ``preserving a little implies preserving
much'' for transitive nep forcings.
\end{list}

\section{Suslin$^+$ and transitive nep forcing}\label{sec:suslinplus}

\subsection*{A Note on Normal $\ZFCx$}

Let us recall the definition of properness:
\begin{Def}\label{def:proper}
$P$ is proper if for some large
regular cardinal $\chi$, for all $p\in P$ and all countable elementary
submodels $N\esm H(\chi)$ containing $p$ and $P$
there is a $q\leq p$ which is
$N$-generic.
\end{Def}

Intuitively, one would like to use elementary submodels of the universe
instead of $H(\chi)$, but for obvious reasons this is not possible.
So one has to show that the properness notion does not depend on the
particular $\chi$ used in the definition, and that essential forcing
constructions are absolute between $V$ and $H(\chi)$ (and $V[G]$ and
$H^{V[G]}(\chi)$). 
So while the choice of $\chi$ is not important, it is not a good idea to
fix a specific $\chi$ (say, $\beth_\omega^+$), since we might for example want
to apply the properness notion to forcings larger than this specific $\chi$.

In Suslin forcing, instead of countable elementary submodels 
arbitrary countable transitive models of some theory $\ZFCx$, so-called
candidates, are used.
Intuitively one would like to use ZFC, but this cannot be done 
for similar reasons. (For example, ZFC does not prove the existence 
of a  ZFC-model.)
%and even if we assume that there are (wellfounded) $ZFC$-models then there will
%always be $ZFC$-models which think that there are no $ZFC$-models.

Again, it turns out that the choice of $\ZFCx$ is of no real importance
(provided it is somewhat reasonable), but we should not 
fix a specific $\ZFCx$.\footnote{We will sometimes require that 
every $\ZFCx$-candidate $M$ thinks that there is a
$\ZFCxx$-candidate $M'$ (and this fails for $\ZFCxx=\ZFCx$), or
that any forcing extension $M[G]$ of $M$ satisfies $\ZFCxx$.}
\newpageforfinallayout
\begin{Def}\label{def:normal}
\begin{itemize}
\item
  $\ZFCminus$ denotes ZFC minus the powerset axiom plus ``$\beth_\omega$ exists''.
\item
  An $\in$-theory $\ZFCx$ is called normal if
  $H(\chi)\vD\ZFCx$ for large regular $\chi$.
\item
  A recursive theory $\ZFCx$ is strongly normal if ZFC proves 
    \[\ex \chi_0\, \fa(\chi>\chi_0\text{ regular })H(\chi)\vD\ZFCx.\]
\end{itemize}
\end{Def} 

We will be interested in strongly normal theories only.
%(Note however that non-normal $\ZFCx$ could be used in the definition of Suslin
%proper and transitive nep as well, which can be useful under specific
%circumstances.) 
Clearly, $\ZFCminus$ is strongly normal. Also,
 if $T$ is strongly normal, then  the theory
$T$ plus ``there is a $T$-candidate'' is strongly normal, and
a finite union of strongly normal theories is strongly
normal.\footnote{This is not true for countable unions, of course:
By reflection, for every finite $T\subset \ZFC$, $\text{Con}(T)$
is strongly normal, but ZFC cannot prove $H(\chi)\vD\text{Con}(\ZFC)$.}

The importance of normality is the following:
If $\ZFCx$ is normal, then forcings that are non-elementary
proper with respect to $\ZFCx$ are proper (see 
the remark after \ref{def:suslinproper}).
However, normal doesn't necessarily mean ``reasonable''. For example, if 
%CH holds in $V$, then $\ZFCx$ plus CH is normal. Also, if 
in $V$ there is no inaccessible, then $\ZFCminus$ plus the negation of the
powerset axiom is normal. 
%(It could even be strongly normal, if one believes
%that inaccessibles are not even {\em consistent}.  To avoid such
%considerations, we could require a strongly normal $\ZFCx$ to be a subset of
%ZFC. However, for this paper this is not relevant.)

As usual, we will (without further mentioning)
assume that certain (finitely many) strongly normal
sentences are in $\ZFCx$. For example, we will state that 
Borel-relations are absolute between candidates and $V$,
which of course assumes that $\ZFCx$ contains enough
of $\ZFCminus$ to guarantee this absoluteness.

%If $Q$ is Suslin and ccc, and in a model $M$ the
%completeness theorem for Keisler-logic
%$\varphi_\text{Keisler}$ is true, then any $Q$-generic filter is 
%generic over $M$ as well (see \ref{lem:cccabs}).
%If $\ZFCx$ is normal, then 
%without loss of generality
% we can assume that $\ZFCx$ contains 
%$\varphi_\text{Keisler}$, and even the following sentence $\varphi_1$:
%``all forcings $R$ such that $2^\card{R}$ exists force that
%$\varphi_\text{Keisler}$ holds.'' (We can assume that 
%$\varphi_1\in\ZFCx$, since 
%it is provable in $\ZFC$ that $\varphi_1$ holds in $H(\chi)$
%for sufficently large regular cardinals $\chi$.
%By the way, we will frequently and without mentioning use the well known fact that for
%large (with respect to $\tau$), regular $\chi$: $H(\chi)\vD p\forc\varphi(\tau)$ iff
%$p\forc H(\chi)^{V[G]}\vD\varphi(\tau)$.)

%Note that for no $\ZFCx$, we can get the folowing 
%``hereditary normality'':
%``Every candidate thinks there exists a ($\ZFCx$-)candidate.''
%(Otherwise we would get an infinite descending chain of candidates.)
%
%As an example for how we will use normality, assume that for a name $\h\in
%H(\al1)$ and all candidates $M$, $M\vD \forc\h\notin V$. Then this is true
%in $V$ as well.  Otherwise, if $V\vD p\forc\h=r$, then for some regular
%$\chi$, $H(\chi)\vD p\forc\h=r$ (since $\h[G]=r$ is absolute
%between $H(\chi)^{V[G]}$ and $V[G]$). Take $N\esm H(\chi)$ countable, $M$ its
%transitive collapse. Then $M$ is a candidate, and $M\vD p\forc\h=r$, a
%contradiction. We will usually abbrevate arguments of this kind by just
%refering to normality.

\subsection*{Candidates, Suslin and Suslin$^+$ forcing}

The following basic setting will apply to all versions 
of Suslin forcings used in this paper (Suslin proper, Suslin ccc,
Suslin$^+$) as well as transitive nep:

We assume that the forcing $Q$ is defined by formulas $\varphi_{\in Q}(x)$ and
$\varphi_\leq(x,y)$, using a real parameter $r_Q$.  Fix a normal $\ZFCx$.  $M$
is called a ``candidate'' if it is a countable transitive $\ZFCx$ model and
$r_Q\in M$.  We denote the evaluation of $\varphi_{\in Q}$ and $\varphi_\leq$
in a candidate $M$ by $Q^M$ and $\leq^M$.  

We further assume that in every candidate $Q^M$ is a
set and $\leq^M$ a partial order on this set; and that $\varphi_{\in Q}$ and
$\varphi_\leq$ are upwards absolute between candidates and $V$.\footnote{This
means that if $M_1$ and $M_2$ are candidates such that $M_1\in M_2$, and if
$q\leq^{M_2} p$, then $q\leq^{M_1} p$ and $q\leq^{V} p$.}

A $q\in Q$ is called $M$-generic (or: $Q$-generic over $M$), if $q\forc \qb G_Q\cap
Q^M\text{ is }Q^M\text{-generic over }M\qe$.

Usually (but not necessarily) 
it will be the case that $p\incomp q$ is absolute between $M$ and $V$.
In this case $q$ is $M$-generic iff 
$q\forc D\cap G_Q\neq \emptyset$
for all $D\in M$ such that
$M\vD\qb D\subseteq Q\text{ dense}\qe$. 
If $p\incomp q$ is not absolute, then this is not enough, since it
does not guarantee that $G_Q\cap Q^M$ is a filter on $Q^M$, i.e.
that it does not contain elements $p$, $q$ such that $M\vD\qb p\incomp q\qe$.
In this case, ``$q$ is $M$-generic'' is equivalent to:
\[q\forc \card{A\cap G_Q}=1\text{ for all }A\in M\text{ such that }M\vD\qb
A\subseteq Q\text{ is a maximal antichain}\qe.\]

We will only be interested in the case $Q\subseteq H(\al1)$.  Assume
$\chi$ is regular and sufficiently large, and $N\esm H(\chi)$ is countable. Let
$i: N\fnto M$ be the transitive collapse of $N$.  Then $i\restriction Q$ is the
identity, and $M$ is a candidate. If $Q$ is proper,  then for every $p\in
Q^M$ there is an $M$-generic $q\leq p$.

Sometimes it would be useful to have generic  conditions for
other candidates (that are not transitive collapses of elementary submodels). 
The first notion of this kind was Suslin proper:

\begin{Def}\label{def:suslinproper}
A (definition of a) forcing $Q$ is Suslin (or: strongly Suslin)
in the parameter $r_Q\in\bR$, if:
%with respect to $\ZFCx$, if:
\begin{enumerate}
\item $r_Q$ codes three $\mS11$ relations, $R^\in_Q$, $R^\leq_Q$ and
	$R^\incomp_Q$.
\item $R^\leq_Q$ is a partial order on $Q=\{x\in{\omega\ho}:\, R^\in_Q(x)\}$
	and $p\incomp_Q q$ iff $R^\incomp_Q(p,q)$.
\setcounter{tmp}{\value{enumi}}	
\end{enumerate}
$Q$ is Suslin proper with respect to some normal $\ZFCx$, if in addition:
\begin{enumerate}\setcounter{enumi}{\value{tmp}} 
\item for every candidate $M$ and every $p\in Q^M$ there is an
    $M$-generic $q\leq p$.
\end{enumerate}
\end{Def}

Remarks:
\begin{itemize}
\item A forcing $Q$ (as a partial order)
	is called Suslin (proper), if there is a
	definition of $Q$ which is Suslin (proper).
\item ``$r_Q$ codes a Suslin forcing'' is a $\mP12$ property.
So if
$Q$ is Suslin in $V$, then $Q$ is Suslin in all candidates 
and all forcing extension of $V$ as well.
In particular, in every candidate $M$, $\leq^M$ is a partial
order on the set $Q^M$ and $p\incomp q$ is equivalent to 
$R^\incomp_Q(p,q)\qe$.\\
However, the formula ``$(\in_Q,\leq_Q,r_Q,\ZFCx)$ codes 
a Suslin proper forcing'' is 
a $\mP13$ statement and in general not absolute.
%(3) will not hold any more in candidates.
%i.e. a Suslin forcing $Q$ that is Suslin proper in $V$ 
%is not necessarily proper in a candidate $M$.
%The statement ``$\leq_Q$ is a (quasi) p.o. on $Q$'' is equivalent to
%$\fa x,y:\, (x\in Q\THEN x\leq x),\ (x\leq y\OR y\leq x\THEN x\in Q),\ 
%(x\leq y\,\AND\,y\leq z\THEN x\leq z) $ which is $\mP12$,
%and ``for every candidate $M$ \dots'' is equivalent to:
%$\fa E\sunseteq \omega \text{ wellfounded such that } (\omega, E)\vD\ZFCx,\,
%\fa p\in \omega$
%which is $Pi^1_2$ as well.
\item If $Q$ is Suslin, then $\incomp$ is a Borel relation, 
and therefore the statement\\
\centerline{$\{q_i:\, i\in\omega\}$ is predense below $p$} 
(i.e.  $p\forc G\cap\{q_i:\, i\in\omega\}\neq \emptyset$)
is $\mP11$ (i.e. relatively $\mP11$ in the $\mS11$ set 
$Q^{(\omega+1)}$).
\item 
If $Q$ is Suslin proper with respect to $\ZFCx$,
and $\ZFCxx$ is stronger than $\ZFCx$, then
$Q$ is Suslin proper with respect to $\ZFCxx$ as well.
\item If $Q$ is Suslin proper, then $Q$ is proper.\\(As mentioned
already, the transitive collapse $M$ of a countable $N\esm H(\chi)$
is a candidate, $Q$ is not changed by the collapse, 
and $q\leq p$ is $M$-generic iff $q\leq p$ is $N$-generic.)
\item The definition of Suslin proper forcing could 
be applied to non-normal $\{\in\}$ theories $\ZFCx$ as well.
This could be useful in other context, but not for
this paper. Obviously such a forcing $Q$ need not be proper any more.
As an extreme example, $\ZFCx$ could contain ``$0=1$''.
Then (3) is immaterial, since there are no candidates, and
every forcing definition $Q$ satisfying (1) and (2) is Suslin proper.
%(but not necessarily proper) with respect to this $\ZFCx$.
\end{itemize}

In \cite{JdSh:292} it is proven that if a forcing $Q$ is Suslin and ccc 
(in short: Suslin ccc), then $Q$ is Suslin proper in a very absolute way:

\begin{Lem}\label{lem:cccabs}
``$Q$ is Suslin ccc'' is a $\mP12$ statement. So in particular, 
if $Q$ is Suslin ccc, then
\begin{enumerate}
\item $Q$ is Suslin ccc in every candidate $M$ and in every forcing extension of $V$.
\item $Q$ is Suslin proper: even $1_Q$ is generic for every candidate.
\end{enumerate}
\end{Lem}

The proof proceeds as follows: Assume $Q$ is Suslin.
Using the completeness theorem $\varphi^\text{Keisler}$ 
for the logic $L_{\om1\omega}(Q)$ (see
\cite{Ke70}) it can be shown \cite[3.14]{JdSh:292}
that ``$Q$ is ccc'' is a Borel statement.
(This requires that $\varphi^\text{Keisler}\in\ZFCx$, which we can
assume since $\varphi^\text{Keisler}$ is strongly normal.)
So if $M$ is a candidate and $M\vD\qb A\subseteq Q\text{ is a maximal
antichain}\qe$, then $M\vD\qb A\text{ is countable}\qe$. And we have already
seen that for $Q$ Suslin and $A$ countable, the statement 
``$A$ is predense'' is $\mP11$ (and therefore absolute). So $A$
is predense in $V$, and $1_Q$ forces that $G_Q$ meets $A$.

Note that (1) and (2) of the lemma are trivially true for a $Q$ that is
definable without parameters (e.g. Cohen, random, amoeba, Hechler), assuming
that $\ZFC\vdash Q\text{ is ccc}$ and $\ZFCx\vdash Q\text{ is ccc}$.

For further reference, we repeat a specific instance of the 
last lemma here:

\begin{Lem}\label{gendownwabs}
If $Q$ is Suslin ccc, $M_1\subseteq M_2$ are candidates, 
and $G$ is $Q$-generic over $M_2$ or over $V$,
then $G$ is $Q$-generic over $M_1$.
\end{Lem}

Cohen, random, Hechler and amoeba forcing are Suslin ccc and Mathias forcing
is Suslin proper.  Miller and Sacks forcing, however, are not, since
incompatibility is not Borel.

This motivated a  generalization of Suslin proper, 
Suslin$^+$ \cite[p. 357]{tools}:
here, we do not require $\incomp$ to be $\mS11$,
so ``$\{q_i:\, i\in\omega\}$ is predense below $p$'' will generally 
not be $\mP11$ any more, just $\mP12$. 
However, we require that there is a 
$\mS12$ relation $\epd$ (``effectively predense'') that 
holds for ``enough'' predense sequences:

\begin{Def}\label{def:suslinplus}
A (definition of a) forcing $Q$ is Suslin$^+$ in the parameter $r_Q$ 
with respect to $\ZFCx$, if:
\begin{enumerate}
\item $r_Q$ codes two $\mS11$ relations, $R^\in_Q$ and $ R^\leq_Q$, and
	an $(\omega+1)$-place $\mS12$ relation $\epd$.
\item In $V$ and every candidate $M$,
	$\leq$ is a partial order on $Q$,
	and $\epd(q_i,p)$ implies ``$\{q_i:\, i\in\omega\}$ 
	is predense below $p$''.
\item for every candidate $M$ and every $p\in Q^M$ there is a
        $q\leq p$ such that
   every dense subset $D\in M$ of $Q^M$ 
   has an enumeration  $\{d_i:\, i\in\omega\}$ 
   such that $\epd(d_i,q)$ holds.
\end{enumerate}
\end{Def}
Again, a partial order $Q$ is called Suslin$^+$ if it has a 
suitable definition.

Clearly, every Suslin proper forcing is Suslin$^+$: $\epd$ can just be
defined by ``$\{q_i:\, i\in\omega\}$ is predense below $p$'', which is even
a conjunction of 
$\mP11$ and $\mS11$, 
 and then the condition \ref{def:suslinplus}(3) is
just a reformulation of 
\ref{def:suslinproper}(3). 

\subsection*{Effective Axiom A}

The usual tree-like forcings are Suslin$^+$.
Here, we consider the following forcings consisting of trees
on $\omega\hko$ ordered by $\subseteq$.
(Usually, Sacks is defined
on $2\hko$, but this is equivalent by a simple density argument.)
For $s,t\in \omega\hko$ we write $s\leq t$ for ``$s$ is an initial segment of
$t$''; for a tree $T\subseteq \omega\hko$ $s\leq_T t$ means $s\leq t$ and
$s,t\in T$; and $s^\frown n$ is the immediate successor of $s$ with last
element $n$.
\begin{itemize}
\item 
Sacks, perfect trees: 
  $(\fa {s\in T})\  (\ex {t\geq_T s})\  (\ex ^{\geq 2}n)\  t^\frown n\in T$.
\item
Miller, superperfect trees: every node has either exactly one or 
  infinitely many immediate successors, and 
  $(\fa {s\in T})\  (\ex {t\geq_T s})\  (\ex^\infty n)\  t^\frown n\in T$.
\item Ros{\l}anowski:  every node has either exactly one or all possible 
  successors, and 
  $(\fa {s\in T})\  (\ex {t\geq_T s})\  (\fa{n\in\omega})\  t^\frown n\in T$.
\item
Laver: let $s$ be the stem of $T$. Then
  $(\fa {t\geq_T s})\  (\ex^\infty n)\  t^\frown n\in T$.
\end{itemize}

In the following, we call Sacks, Miller and Ros{\l}anowski ``Miller-like''.
Clearly, ``$p\in Q$'' and ``$q \leq p$'' are Borel (but $p\incomp q$ is
not).\footnote{Alternatively, $Q$ could of course be defined as the set of
trees just {\em containing} a corresponding set, then $x\in Q$ is $\mS11$, and
for the Miller-like forcings two compatible elements $p$, $q$ have a canonical
lower bound, $p\cap q$.}

For Sacks, there is a proof of the Suslin$^+$ property in \cite{tools}
and \cite{pmbc}
using games. 
However, in the same way as the ``canonical'' proof of properness of these
forcings uses Axiom A, the most transparent way to prove
Suslin$^+$ uses an effective version of Axiom A:

Baumgartner's Axiom A \cite{Baumgartner}  for a forcing $(Q,\leq)$
can be formulated as follows:
There are relations $\leq_n$ such that 
\begin{enumerate}
\item $\leq_{n+1}\,\subseteq\,\leq_n\,\subseteq\,\leq$.
\item Fusion: if $(a_n)_{n\in\omega}$ is a sequence of elements of $Q$
  such that $a_{n+1}\leq_n a_n$ then there is an $a_\omega$ such that 
  $a_\omega \leq a_n$ for all $n$.
\item If $p\in Q$, $n\in \omega$ and $D\subseteq Q$ is dense
then there is a $q\leq_n p$ and a countable subset $B$ of $D$
which is predense under $q$.
\end{enumerate}
Remarks:
\begin{itemize}
\item
  Actually, this is a weak version of Axiom A, usually 
  something like $a_\omega \leq_n a_n$ will hold in (2). 
\item  
  It is easy to see that in (3), instead of 
  ``and $D\subseteq Q$ is dense'' 
 we can equivalently use ``and $D\subseteq Q$ is open dense''
 (or maximal antichain).
\end{itemize}

Now for ``effective Axiom A'' it is required that the $B\subseteq D$ in (3)
is {\em effectively} predense below $q$, not just predense. 
Then Suslin$^+$ follows.
To be more exact:

\begin{Def}\label{def:effectiveaxioma}
$Q$ satisfies effective Axiom A (in the parameter $r_Q$ 
with respect to $\ZFCx$), if
\begin{enumerate}
\item $r_Q$ codes $\mS11$ relations, $R^\in_Q$, $R^\leq_Q$, and
	$\mS12$ relations $\leq_Q^n$ ($n\in\omega$) and 
	an $(\omega+1)$-place $\mS12$ relation $\epd$.
\item In $V$ and every candidate $M$,
	$\leq$ is a partial order on $Q$ %=\{x\in{\omega\ho}:\, R^\in_Q(x)\}$
	and $\epd(q_i,p)$ implies that $\{q_i:\, i\in\omega\}$ is predense below $p$.
\item Fusion: %In $V$, 
     For all $(a_n)_{n\in\omega}$ %\in{Q\ho}$ 
     such that
     $a_{n+1}\leq_n a_n$ there is an $a_\omega$ such that $a_\omega \leq a_n$.
\item In all candidates, if $p\in Q$, $n\in\omega$ and $D\subseteq Q$ is dense then
   there is  a $q\leq_n p$ and a sequence $(b_i)_{i\in\omega}$ of
   elements of $D$ such that $\epd(b_i,q)$ holds.
\end{enumerate}
\end{Def}
Again, a partial  order  $Q$ satisfies effective Axiom A
if it has  a suitable definition.

\begin{Lem}
If the partial order $Q$ satisfies effective Axiom A, 
then $Q$ is Suslin$^+$.
\end{Lem}

\begin{proof}
First we define $\epd'(p'_i,q')$ by
\[(\ex {q\geq q'})\, (\ex {\{p_i\}\subseteq \{p'_i\}})\, \epd(p_i,q).\]
Clearly, this is a $\mS12$ relation coded by $r_Q$ 
satisfying \ref{def:suslinplus}(2).
Let $M$ be a candidate, and let $\{D_i:\, i\in\omega\}$ list the dense sets
of $Q^M$ that are in $M$. Pick an arbitrary $a_0=p\in Q^M$.
We have to find a
$q\leq p$ satisfying \ref{def:suslinplus}(3) with respect to $\epd'$. 
Assume we have already constructed $a_n$.
In $M$, according to (4) using $D_n$ as $D$, we find
an $a_{n+1}\leq_n a_n$ and $\{b^n_i:\, i\in\omega\}\subseteq D_n$ 
such that $\epd(b^n_i,a_{n+1})$ holds (in $M$ and 
therefore by absoluteness in $V$).
In $V$ pick $q=a_\omega$ according to (3).
\end{proof}

The usual proofs that the forcings defined above  satisfy Axiom A
also show that they satisfy the effective version. 
To be more explicit: Let $Q$ be one of the  forcings.
We define (for $p$, $q\in Q$, $n\in\omega$):
\begin{itemize}
\item $\spl(p)=\{s\in p:\, (\ex^{\geq 2}n\in\omega)\ s^\frown n\in p\}$.
\item $\spl(p,n)=\{s\in \spl(p):\, (\ex^{=n}{t\leq s})\, t\in\spl(p)\}$.\\
	(So $s\in\spl(p,n)$ means that $s$ is the $n$-th splitting node
	along the branch $\{t\leq s\}$. In particular, $\spl(p,0)$
	is the singleton containing the stem of $p$.)
\item $q\leq_n p$, if $q\leq p$ and $\spl(q,n)=\spl(p,n)$.\\
	(So $q\leq_0 p$ if $q\leq p$ and $q$ has the same stem as $p$.)
\item For $s\in p$, $p^{[s]}=\{t\in p:\, {t\leq s}\OR{s\leq t}\}$.
\item $F\subseteq p$ is a front (or: $F$ is a front in $p$), 
	if it is an antichain meeting every branch of $p$.
\item $\epd(q_i,p)$ is defined by: There is a front $F\subseteq p$
  such that $\fa{t\in F}\, \ex{i \in\omega}:\, q_i=p^{[t]}$.
\item For Miller-like forcings, effectively predense could also be
 defined as $\epd'(q_i,p)\ \DEFIFF\ \ex n\, \fa{s\in\spl(p,n)}\, \ex i:\, 
 q_i=p^{[s]}$.
\end{itemize}

Clearly, $\spl(p)$, $\spl(p,n)$, $p^{[s]}$ and $\epd'$ are Borel,
``$F$ is a front'' is $\mP11$, therefore $\epd$ is $\mS12$.
The following facts are easy to check ($p,q \in Q$):
\begin{itemize}
\item If $s\in p$, then $p^{[s]}\in Q$. 
\item If $F\subset p$ is a front and $q\comp p$, then 
	$q\comp p^{[s]}$ for some $s\in F$.
\item $\spl(p,n)$ is a front in $p$.
\item For $(q_n)_{n\in\omega}$ such that $q_{n+1}\leq_n q_n$, there is a canonical
	limit $q_\omega\in Q$ and $q_\omega\leq_n q_n$.
%\end{itemize}
%\setcounter{tmp}{\value{enumi}}	
%\begin{itemize}%\setcounter{enumi}{\value{tmp}} 
\item 
  If $Q$ is Miller-like, and
  if $F\subset p$ is a front, and $\fa{s\in F}, p_s\in Q,\, p_s\subseteq 
	p^{[s]}$, then $\bigcup_{s\in F} p_s\in Q$, and
	$\bigcup_{s\in F}p_s\subseteq p$.
\item 
  If $Q$ is Laver, and if $F\subset p$ is a front, and 
  $\fa{s\in F}, p_s\in Q$ has stem $s$, then
  %and the stem of $p_s$ is $s$, then
  $\bigcup_{s\in F} p_s\in Q$, and
        $\bigcup_{s\in F}p_s\subseteq p$.
\end{itemize}  

$\leq_n$ and $\epd$ defined as above satisfy the requirements
\ref{def:effectiveaxioma} for effective Axiom A:\\
(1)--(3) are clear.\\
For Miller-like forcings, (4) is proven as follows:
Assume $D\subseteq Q$ is dense and $p\in Q$. For all $s\in \spl(p,n+1)$,
$p^{[s]}\in Q$, so there is a $q^s\subseteq p^{[s]}$ such that $q^s\in D$.
Now set \mbox{$q:=\bigcup_{s\in F} q^s\in Q$}. Then $q\leq_n p$, and the set
$\{q^s:\, s\in F\}\subseteq D$ is effectively predense below $q$ according
to the definition of $\epd'$ (or $\epd$).\\
For Laver, we have to define a rank of nodes: Assume $D$ is dense,
and $p_0$ a condition with stem $s_0$, $s\geq s_0$, and $s\in p_0$. We 
define $\rank_D(p_0,s)$ as follows:\\
\phantom{xxxx}If there is a $q\subseteq p_0$ such that $q\in D$ and 
    $q$ has stem $s$, then $\rank_D(p_0,s)=0$.\\
\phantom{xxxx}Otherwise $\rank_D(p_0,s)$ is the minimal $\alpha$ such that
for infinitely many immediate\\
\phantom{xxxx}successors $t$ of $s$ the following holds:
    $t \in p_0$ and $\rank_D(p_0,t)<\alpha$.\\
$\rank_D$ is well-defined for all nodes $\geq s_o$ in $p_0$:\\
Assume towards a contradiction that $\rank_D(p_0,s)$ is undefined.
Then
\[q:=\{s'\in p_0^{[s]}:\, s'\leq s\text{ or }\rank_D(p_0,s')\text{ undefined}\}\]
is a Laver condition stronger than $p_0$. Pick a 
$q'\leq q$ such that $q'\in D$. Let $s'$ be the stem of $q'$.
Then $\rank_D(p,s')=0$, $s'\geq s$ and $s'\in q$,  a contradiction.\\
Now define $q'\leq p_0$ inductively. First add all
$s\leq s_0$ to $q'$.  Assume $s\in q'$ and $s\geq s_0$.
Then we add infinitely many immediate successors $t\in p_0$ of $s$ to $q'$.
If $\rank_D(p,s)\neq0$, we additionally require that
$\rank_D(p,t)<\rank_D(p,s)$ for each of these $t$
(this is possible by the
definition of $\rank_D(p,s)$). So the $q'$ constructed this way
is a Laver condition with the same stem $s_0$ as $p_0$. Also, along
every branch of $q'$, $\rank_D(p,s)$ is strictly decreasing (until
it gets 0), therefore there is a front $F_0$ in $q'$ such that for all 
$s\in F_0$, $\rank_D(p,s)=0$. That means that for all $s\in F_0$ there is a 
$q^s\leq p_0$
such that $q^s\in D$ and $q^s$ has stem $s$. 
Define $q_0$ to be $\bigcup_{s\in F_0} q^s$.
Clearly $q_0\leq p_0$, $q_0$ has the same stem $s_0$ as $p_0$, 
$F_0$ is a front in $q_0$ and for every $s\in F_0$, $q_0^{[s]}\in D$.\\
Given a Laver condition $p$ and $n\in\omega$, define for every
$p_0\in\spl(p,n)$ a $q_0$ as above, and let $q$ be the union of these
$q_0$, and $F$ the union of the according $F_0$. 
Then $q\leq_n p$, and for every $s$ in the front $F\subset q$,
$q^{[s]}\in D$.
This finishes the proof of effective Axiom A for Laver.
%the value of $\rank_D$ is strictly decreasing along branches,
%therefore $F=\{s\in p:\, \rank_D(p,s)=0\}$ is a front.
%So for each $s\in F$ there is a $p_s\leq p$ in $D$ with stem $s$.
%So $\{p_s:\, s\in F\}$ is effectively predense below 
%$\bigcup_{s\in F}p_s$. % which is $\leq p$

It is clear that the same proof of effective Axiom A works for other tree
forcings as well, for example for all finite-splitting lim-sup tree forcings.
(In \cite[1.3.5]{RoSh:470} such forcings are called $\mathbb{Q}^\text{tree}_0$.)

\subsection*{Transitive nep}

So we have seen that Suslin ccc implies Suslin
proper, which implies Suslin$^+$. 
For the proof of
the main theorem \ref{maintheorem}, even less than Suslin$^+$ is 
required:\footnote{Actually, for the main theorem even less than
nep would be sufficient: we need generic conditions for candidates $M$ that are
internal set forcing extensions of transitive collapses of elementary submodels
only (not for all candidates). However, this restriction doesn't seem to
lead to a natural nep notion.}
A forcing definition $Q$ (using the parameter $r_Q$) is transitive nep
(non-elementary proper),
if 
\begin{itemize}
\item ``$p\in Q$'' and ``$q\leq p$'' are upwards absolute between candidates 
and $V$.
\item In $V$ and all candidates, $Q\subseteq H(\al1)$
and ``$p\in Q$'' and ``$q\leq p$'' are absolute
between the universe and $H(\chi)$ (for large regular $\chi$).
\item For all candidates $M$ and $p\in Q^M$
there is a $q\leq p$ forcing that $G\cap Q^M$ is $Q^M$-generic over $M$.
\end{itemize}

Recall our initial consideration:
In proper forcing, we get the properness condition for (collapses of)
elementary submodels only, but we would like to have it for non-elementary
models as well. (This is the reason for the name ``non-elementary proper''.) So
transitive nep captures this consideration with little additional assumptions.

There is also a (technically more complicated) version of nep for
non-elementary and non-transitive candidates, defined in \cite{Sh:630}, which
makes it possible for long iterations to be nep (transitive nep requires
$Q\subseteq H(\al1)$).  The main theorem \ref{maintheorem} of this paper 
holds for this general notion of nep as well (with nearly the same proof).

For every countable transitive model, 
$M\vD\qb p\forc \varphi(\tau)\qe$ iff for all $M$-generic $G$ containing $p$,
$M\vG\vD\qb\varphi(\tau\vG)\qe$.
If $Q$ is nep and $M$ a candidate, then
$M\vD\qb p\forc \varphi(\tau)\qe$ iff for all $M$- {\em and} $V$-generic $G$
containing $p$,
$M\vG\vD\qb\varphi(\tau\vG)\qe$:\\
One direction is clear. For the other, assume $M\vD\qb p'\leq p,\, p'\forc
\lnot \varphi(\tau)\qe$.  Let $q\leq p'$ be $M$ generic.  Then for any
$V$-generic $G$ containing $q$, $G$  is $M$-generic as well and
$M\vG\vD\qb\lnot\varphi(\tau\vG)\qe$.

We will need the following instance of Shoenfield-Levy absoluteness:
\begin{Lem}\label{lem:shoenfield}
Let $x\in H(\al1)$.
Then ``there is a candidate $M$ containing $x$ such that
$M\vD\varphi(x)$'' is $\mS12$
%Let $V_1\subseteq V_2$ be two transitive models of $\ZFC$, $\om1\subset
%V_1$, $V_1\vD x\in H(\al1)$.
(and therefore absolute between universes with the same $\om1$).
%is absolute between $V_1$ and $V_2$.
\end{Lem}

%This is just $\Sigma_1$ (Shoenfield-Levy) absoluteness. 

All in all we get the following implications:
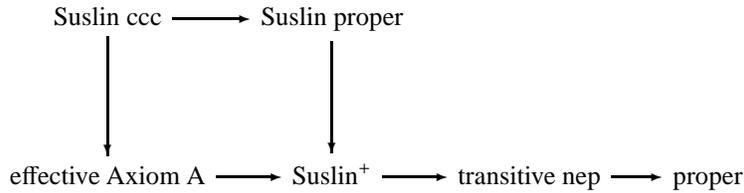
\begin{diagram}
\text{Suslin ccc}         & \rTo & \text{Suslin proper}  \\
     \dTo                 &      &      \dTo             \\
\text{effective Axiom A}  & \rTo &  \text{Suslin}^+    & \rTo & \text{transitive nep} &\rTo&\text{proper}\\
\end{diagram}
\newpageforfinallayout
\section{Suslin ccc ideals}\label{sec:suslinccc}

The set of Borel codes (or Borel definitions) will be denoted by
``$\BC$''. So $\BC$ is a set of reals.
For $A\in \BC$ we denote 
the set of reals that satisfy the definition $A$
(in the universe $V$) with  $A^V$.

If $Q\subseteq H(\al1)$ is ccc, 
then a name $\n\tau$ for an element of $\omega\ho$
can be transformed into an equivalent 
hereditarily countable name $\h$:
for every $n$, pick a maximal antichain $A_n$ deciding $\n\tau(n)$,
then $\h\DEFEQ\{(p,(n,m)):\, p\in A_n, p\forc\n\tau(n)=m\}$ is equivalent to $\n\tau$.

From now on, we will assume the following:

\begin{Asm}\label{basicasm}
$Q$ is a Suslin ccc forcing, 
$\h$ is a hereditarily countable name
coded by $r_Q$,
$\forc_Q \h\in{\omega\ho}\setminus V$, and in all candidates:
$\{\llbracket \h(n)=m\rrbracket:\, n,m\in\omega\}$ generates $\ro(Q)$.
\end{Asm}

``$X$ generates $\ro(Q)$'' means that there is no proper
sub-Boolean-algebra $B\supseteq X$ of $\ro(Q)$ such that 
$\sup_{\ro(Q)}(Y)\in B$ for all $Y\subseteq B$.

\begin{Lem}\label{lem:vvstrich}
	This assumption is a $\mP12$ statement.
\end{Lem}
%(Here we assume that $\ZFC0\in M$, otherwise ``for all candidates'' cannot
%be formulated in $M$.)
\begin{proof}
``$Q$ is Suslin ccc'' is $\mP12$ according to \ref{lem:cccabs}.
For $x\in H(\al1)$, a statement of the form 
``every candidate thinks $\varphi(x)$''
is $\mP12$ (cf. \ref{lem:shoenfield}).
$\forc_Q (\h\in{\omega\ho}\setminus V)$ holds in $V$ iff it holds in 
every candidate:
  If $M\vD p\forc\h=r$, 
  then this holds in $V$ as well: For Suslin ccc forcings,
  every $V$-generic filter is $M$-generic, and
  $\h=r$ is absolute. The other direction follows 
  from normality.
\end{proof}

\begin{Lem}\label{lem:Iqisabs}
	For $A\in\BC$, ``$q\forc \h\in A^{V[G_Q]}$'' is $\mD12$.
\end{Lem}

Remark: \cite[2.7]{BaBo} gives a general result for $\mS1n$ formulas.

\begin{proof}
  For any candidate $M$ containing $q$ and $A$,  ``$q\forc \h\in A$''
  is absolute between $V$ and $M$:
  If $q\in G$ is $V$-generic, then
  it is $M$-generic as well (since $Q$ is Suslin ccc), and
  $\h[G]\in A$ is absolute between $M[G]$ and $V[G]$.

  So $q\forc \h\in A$ iff for all candidates $M$, 
  $M\vD q\forc \h\in A$ (a $\mP12$ statement) iff for some candidate
  $M$: $M\vD q\forc \h\in A$ (a $\mS12$ statement). 
\end{proof}

\begin{Lem}\label{uniquniq} The statement\\
\centerline{$\{\llbracket \h(n)=m\rrbracket:\, n,m\in\omega\}$ generates
$\ro(Q)$} holds in $M$ iff the following holds (in $V$):\\
\centerline{if $G_1,G_2\in V$ are $Q$-generic over $M$ and $G_1\cap M\neq G_2\cap M$, then
$\h[G_1]\neq \h[G_2]$.}
\end{Lem}

\begin{proof}
If $\{\llbracket \h(n)=m\rrbracket:\, n,m\in\omega\}$ generates $\ro(Q)$, then
$G\cap Q\rzM$ can be calculated (in $M\vG$) from $\h\vG$.
On the other hand, let (in $M$) 
$B=\ro(Q)$, $C$ the proper complete sub-algebra
generated by $\llbracket \h(n)=m\rrbracket$.
Take $b_0\in B$ such that no $b'\leq b_0$ is in $C$, and set
\[c=\inf\{ c'\in C:\, c'\geq b_0 \},\quad b_1=c\setminus b_0.\]
So for all $c'\in C$, $c'\comp b_0$ iff $c'\comp b_1$.
Let $G_0$ be $B$-generic over $M$ such that $b_0$ in $G$. Then
$H=G_0\cap C$ is $C$-generic. 
In $M\vH$, $b_1\in B/H$. So there is a  $G_1\supset H$ containing 
$b_1$.
\end{proof}

%Remark: So for a certain $M$, the statement is $\mP11$. Therefore 
%this part of asummtion \ref{} is $\mP12$, another proof of absoluteness

\enlargethispage{1cm}
\begin{Def} The Suslin ccc ideal corresponding to $(Q,\n\eta)$:
\begin{itemize}
\item 
$\IB=\{A\in \BC:\, \forc_Q\h\notin A^{V[G_Q]}\}$.
\item
%A Borel set $A$ is in is in $I$, if
%$\forc_Q\h\notin A$ (where $A$ is interpreted as a Borel code, 
%not as a set).
%\item
%$X\subseteq {\ho \omega}$ is in $I$, if
%there is a Borel $A\in I$ such that $X\subseteq A$.
$I=\{X\subseteq {\omega\ho}:\, \ex A\in \IB:\, A^V\supseteq X\}$.
\item 
$X\in \Ip$ (or: $X$ is positive) means $X\notin I$, and 
$X$ is of measure 1 means ${\omega\ho}\setminus X\in I$.\\
$\IpB:=\BC\setminus \IB$.
\end{itemize}
\end{Def}

Note that we use the phrases ``of measure 1'', ``null'' and
``positive'' for all Suslin ccc ideals, not just for the 
Lebesgue null ideal.
For example, if
%if $\bB$ is the random algebra
%then $I_\bB$ are the null sets, $\Ip_\bB$
%are the non-null sets and $\coI_\bB$ means
%measure 1.
%If 
$\bC$ is Cohen forcing,
then the null sets are the meager sets,
and a set has ``measure 1'' if it is co-meager.

Clearly $A\in \IB$ iff $A^V\in I$.

An immediate consequence of lemma \ref{lem:Iqisabs} is
\begin{Cor}
For $A\in\BC$,  ``$A\in \IB$'' is $\mD12$.
\end{Cor}

So for Borel sets, being null is absolute.

\begin{Lem}\label{Phi}
$I$ is a $\sigma$-complete ccc ideal containing all singletons, and
there is a surjective $\sigma$-Boolean-algebra homomorphism
$\phi:\text{Borel}\fnto \ro(Q)$ with kernel $I$, i.e.
$\ro(Q)$ is isomorphic to $\text{Borel}/I$ as a complete Boolean algebra.
\end{Lem}

ccc means: there is no uncountable family $\{A_i\}$ such that 
$A_i\in \Ip$ and $A_i\cap A_j\in I$ for $i\neq j$
(or equivalently: $A_i\cap A_j=\emptyset$).

\begin{proof}
$\sigma$-complete is clear: 
If $X_i\subseteq A_i\in I$, and $\forc \h\notin A_i$ for all $i\in\omega$,
then $\forc \h\notin \bigcup A_i\supseteq \bigcup X_i$.\\
%$I$ is ccc follows form the isomorphism of $\ro(Q)$ and $\text{Borel}/I$.
%$\ro(Q)$ is isomorphic to $\text{Borel}/I$ (that will follow from the proof
%of \ref{jglo}: ) 
%and $Q$ is ccc,so $I$ is ccc.
For $A\in  \BC$, define $\phi(A)=\llbracket \h\in A^{V[G]}\rrbracket_{\ro(Q)}$.
Then \mbox{$\phi({\omega\ho}\setminus A)=\lnot \phi(A)$},
\mbox{$\phi(\U A_i)=\sup\{\phi(A_i)\}$}, and
if $A\subseteq B$, then $\phi(A)\leq \phi(B)$.
If $\phi(A) \leq \phi(B)$, then $\forc \h\notin (A\setminus B)$,
so $A\setminus B\in I$.
Since $\h$ generates $\ro(Q)$ (in all candidates, and
therefore in $V$ as well by normality)
and since $Q$ is ccc, $\ro(Q)=\phi''\text{Borel}$.
So $\phi:\text{Borel}\fnto \ro(Q)$ 
is a surjective $\sigma$-Boolean-algebra homomorphism.
The kernel is the $\sigma$-closed ideal $I$, so $\text{Borel}/I$ is 
isomorphic to $\ro(Q)$ as a $\sigma$-Boolean-algebra, and
(since $\ro(Q)$ is ccc), even as complete Boolean algebra.
\end{proof}

\begin{Def}
%For $q\in Q^M$, $\es$  is called $(Q,\h)$-generic over $M$ containing $q$
%($\es\in\Gen(M,q)$), if there is a $G\in V$ $M$-generic such that $q\in G$ and
%$\h[G]=\es$. $\Gen(M)=\Gen(M,1_Q)$.
$\es$  is called generic over $M$ ($\es\in\Gen(M)$),
if there is an 
$M$-generic 
$G\in V$ 
such that $\h[G]=\es$.
\end{Def}

According to \ref{uniquniq}, this $G$ is unique (on $Q\cap M$).
For example, if $Q$ is random, then $\Gen(M)$ is the set of random 
reals over $M$.

%``$\llbracket \h\in B\rrbracket=q$'' %($B$ Borel code) is $\mD12$.
$\llbracket \h\in B\rrbracket=q$ is equivalent to
\[q\forc \h\in B \text{ and if }p\incomp q\text{ then } p\forc \h\notin B,\]
which is $\mP12$ (because of lemma \ref{lem:Iqisabs} and the fact that
$p\incomp q$ is Borel).
For $q\in Q$ we denote a $B$ such that $\llbracket \h\in B\rrbracket=q$ by
$B_q$. Of course $B_q$ is not unique, just unique modulo $I$.
$q \forc \h \in A$  iff
$\forc ( \h \in B_q \THEN \h\in A)$, i.e. iff
$\forc \h\notin B_q\setminus A$. So we get
%\begin{Lem}\label{hilfslemma3}
$q\forc \h\notin A$ iff $A\cap B_q\in I$, and 
$q\forc \h\in A$ iff $B_q\setminus A\in I$.
%\end{Cor}

If $M$ is a candidate, then because of lemma \ref{lem:vvstrich}
the assumption \ref{basicasm} holds in $M$, 
so $M$ knows about the isomorphism $\ro(Q)\fnto \text{Borel}/I$
and in $M$ there is a $B^M_q$ as above.

\begin{Lem}\label{lem:geneqavoidsnull}
Let $M$ be a candidate and $q\in Q\cap M$. Then
\begin{enumerate}
\item
$\Gen(M)= {\omega\ho}\setminus\U\{A^V:\, A \in \IB\cap M\}$.
\item 
$\{\h\vG:\, G\in V\text{ is } M\text{-generic and }  q\in G\}=$\\
$={\omega\ho}\setminus \U\{A^V:\, A\in\BC\cap M, q\forc \h\notin A^{V[G_Q]}\}
=\Gen(M)\cap B^M_q$.
\item 
$\Gen(M)$ is a Borel set of measure 1.
\end{enumerate}
\end{Lem}

For example, if $Q$ is random forcing, this just says that 
$\es$ is generic (i.e. random) over $M$ iff for all Borel codes 
$A\in M$ of null sets, $\es\notin A^V$. 

\begin{proof}
(1) is just a special case of (2).

(2) Set 
\begin{eqnarray*}
X&:=& {\omega\ho}\setminus \U\{A^V:\, A\in\BC\cap M, q\forc \h\notin A^{V[G_Q]}\},\text{ and}\\
Y&:=& \{\h\vG:\, G\in V\text{ is } M\text{-generic and } q\in G\}.
\end{eqnarray*}
Assume $\es\in Y$. Let $G$ be $M$-generic such that $q\in G$ and
$\h\vG=\es$. If $M\vD q\forc \h\notin A^{V[G_Q]}$, 
then $M\vG\vD \es\notin A^{M\vG}$,
i.e. $\es\notin A^V$. So $\es \in X$.

If $\es \in X$, use
(in $M$) the mapping 
$\phi:\text{Borel}\fnto\ro(Q)$ ($A\mapsto \llbracket \h\in A\rrbracket$).
If $\phi(A) \leq \phi(B)$, then $\forc \h\notin
(A\setminus B)$, so by our assumption, $\es \notin (A\setminus B)$.  Given
$\es$, define $G$ by $\phi(A)\in G$ iff $\es\in A$.  $G$ is well defined:
If $\es\in A\setminus B$, then $\phi(A) \neq \phi(B)$.  We have to show that
$G$ is a generic filter over $M$:
If $\phi(A_1),\phi(A_2)\in G$, then $\es\in A_1\cap A_2$,
	so $\phi(A_1)\wedge \phi(A_2)\in G$. 
If $\phi(A)\leq \phi(B)$, then $\es \notin (A\setminus B)$, 
	so $\phi(A)\in G\THEN \phi(B)\in G$.
Since $\phi(\emptyset)=0$, and $\es\notin \emptyset$, $0\notin G$.
If $\sup(\phi(A_i))\in G$, $(A_i)\in M$, then $\es\in \U A_i$,
i.e. for some $i$, $\phi(A_i)\in G$.
Since $q\forc \h \notin {\omega\ho}\setminus B^M_q$,
$\es\notin {\omega\ho}\setminus B^M_q$, i.e. $\es\in B^M_q$,
and since $\phi(B^M_q)=q$, $q\in G$, so $\es\in Y$.
So we have seen that $Y=X\subseteq \Gen(M)\cap B^M_q$.

If $\es\in \Gen(M)\cap B^M_q$, witnessed by  $G$,
then $\h[G]\in B^M_q$,  so $q\in G$
(since $q=\llbracket \h\in B^M_q\rrbracket$), 
i.e. $\es\in Y$.

(3) follows from 1, since $I$ is $\sigma$-complete.
\end{proof}

Note that if $Q$ is not ccc, then our definition of $I$ does not 
lead to anything useful. For example, if $Q$ is Sacks forcing,
then $I_Q$ is the ideal of countable sets, and clearly 
lemma \ref{lem:geneqavoidsnull} does not hold any more.
There are a few possible definitions for
ideals  generated by non-ccc forcings, see for example
\cite{tomekbook}. For tree-forcings $Q$, a popular
ideal is the following: A set of reals $X$ is in $I$, if for 
every $T\in Q$ there is a $S\leq_Q T$ such that $\lim(S)\cap X=\emptyset$.
In the case of Sacks forcing this ideal is called Marczewski ideal, it 
is not ccc, and a Borel set $A$ is in $I$ iff $A$ is countable.

%%%%%%%%%%%%%%%%%%%%%%%%%%%%%%%%%%%%%%%%%%%%%%%%%%%%%%%%%%%%%%
% PRESERVATION

\section{Preservation}\label{sec:preservation}

\begin{Def}
\begin{itemize}
\item
    $P$ is Borel $I^+$-preserving, if for all $A\in \IpB$,
    $\forc_P A^V\in \Ip$. 
\item
    $P$ is $I^+$-preserving, if for all $X\in\Ip$,
    $\forc_P \std{X}\in \Ip$.
\end{itemize}
\end{Def}

For example, if $Q$=random, then random forcing is $I^+$-preserving, and Cohen
forcing is not Borel $I^+$-preserving. If $Q$=Cohen, then Cohen forcing is
$I^+$-preserving, and  random forcing is not Borel $I^+$-preserving. 

Note that being Borel $I^+$-preserving is stronger than just 
``$\forc_P V\cap {\omega\ho}\notin I$''. For example, 
set $X:=\{x\in{\omega\ho}:\, x(0)=0\}$ and  $Y:={\omega\ho}\setminus X$.
Let $Q$ be the forcing that adds a real $\h$ such that 
$\h$ is random if $\h\in X$ and $\h$ is Cohen otherwise.
Clearly, $Q$ is Suslin ccc.
$A\in I$ iff ($A\cap X$ is null and $A\cap Y$ is meager).
So if $P$ is random forcing, then $\forc_P({\omega\ho}^V \notin I\ \AND\ Y^V\in
I)$. Note that in this case %, for any candidate $M$, 
a $Q$-generic real $\es$ over $M$ will still be generic after
forcing with $P$ if $\es\in X$, but not if $\es \in Y$.

However, if $P$ is homogeneous in a certain way with respect to $Q$, then Borel
$I^+$-preserving and ``$\forc_P V\cap {\omega\ho}\notin I$'' are equivalent
(see \cite{Sh:630} or \cite[3.2]{KrSh:828} for more details).

Also, Borel $I^+$-preserving and $I^+$-preserving are generally not equivalent,
not even if $P$ is ccc. The standard example is the following: Let
$Q$ be $\bC$ (i.e. Cohen forcing, so $I$ is the ideal of meager sets). 
We will construct a 
forcing extension $V'$ of $V$ and a ccc forcing $P\in V'$ such that
$P$ is Borel $I^+$-preserving but not $I^+$-preserving (in $V'$):

Let $\bC_{\om1}$ be the forcing adding $\al1$ many Cohen reals
$(c_i)_{i\in\om1}$, i.e. $\bC_{\om1}$ is the set of all finite partial 
functions from $\omega\times\om1$ to $2$.
Then in any $\bC_\om1$-extension $V[(c_i)_{i\in\om1}]$
%:  
%\begin{itemize}
%\item 
the Cohen reals $\{c_i:\, i\in\om1\}$ are a Luzin set\footnote{$C$
is a Luzin set if $C$ is uncountable and the intersection of $C$ 
with any meager set is countable.}
and for all 
non-meager Borel sets $A$, $A\cap \{c_i:\, i\in\om1\}$ is uncountable. 
If $r$ is random over $V$, and $(c_i)_{i\in\om1}$ 
is $\bC_\om1$-generic over $V[r]$, then
$(c_i)_{i\in\om1}$ is $\bC_\om1$-generic over $V$ as well. 
So the ccc forcing $\bB \ast \bC_\om1$%^{V[G_C]}$ 
can be factored as 
$\bC_\om1\ast \n{P}$, where $\n{P}$ is (a name for a) ccc forcing. 
Set $V':=V[(c_i)_{i\in\om1}]$ and  $V''=V'[G_P]=V[r][(c_i)_{i\in\om1}]$. 
Then in $V'$,
$P=\n{P}[(c_i)_{i\in\om1}]$ is ccc
and Borel $I^+$-preserving, 
${\omega\ho}\cap V \notin I$,
but $P\forc {\omega\ho}\cap V\in I$.

\begin{Def}\label{def:presgen}
\begin{itemize}
\item
    For $p\in P^M$, 
    $\es$  is called absolutely $(Q,\h)$-generic with respect to $p$
    ($\es\in\aGen(M,p)$), if
    there is an $M$-generic 
    $p'\leq p$ forcing that
    $\es\in\Gen(M\vG)$.
\item
    $P$ preserves generics for $M$ if for all $p\in P^M$,
    $\Gen(M)=\aGen(M,p)$. (I.e.
    every $M$-generic real could still be $M\vG$-generic in an extension.)
\end{itemize}
\end{Def}

Note that $\aGen(M,p)\subseteq\Gen(M)$ by \ref{gendownwabs} 
(or \ref{lem:geneqavoidsnull}).

\begin{Lem}
If $P$ preserves generics for (the transitive collapse of) 
unboundedly many countable $N\esm H(\chi)$, 
then $P$ is $I^+$-preserving.
\end{Lem}

Here, unboundedly many means that for all countable $X\subset {\omega\ho}$ 
there is an $N\esm H(\chi)$ countable containing $X$ and $P$ 
with the required property.

Remark: The lemma still holds if $Q$ is any ccc forcing (i.e.
not Suslin ccc. Then
$N$ is not collapsed but used directly as in usual proper forcing
theory).

\begin{proof}
Assume $p\forc_P X\subseteq \n{A}\vGP\in I$, i.e.
$p\forc_P \forc_Q \h \notin \n{A}\vGP^{V[G_P][G_Q]}$. Let $N\esm H(\chi)$
contain $P,X,\n{A},Q,p$. Let $M$ be the collapse of $N$ and 
$\es\in\Gen(M)$,
$p'\leq p$ $M$-generic
such that  $p'\forc\es\in\Gen(M\vGP)$. Let
$G$ be $V$-generic,
$p'\in G$. 

Then $V[G]\vD M\vGP\vGQ\vD \es\notin A\supseteq X$, so 
$V\vD\es\notin X$. Therefore $\Gen(M)\cap X=\emptyset$.
$\Gen(M)$ is of measure 1, therefore $V\vD X\in I$.
\end{proof}

%Before we can formulate the main result we have to clarify the
%degree of normality we need:
%\begin{Def}\label{def:Pnormal}
%$Q$ is called $P$-normal, if $\ZFCx$ is normal
%and $\forc_P (\ZFCx\text{ normal})$, i.e. 
%for large regular $\chi$, $\varphi\in\ZFCx$: 
%$H(\chi)\vD\qb \forc_P\varphi\qe$.
%\end{Def}

%Recall that we assume that $\ZFCP$ is normal
%(so 
%without loss of generality
% it contains the sentence $\varphi_1\AND \varphi_\text{Keisler}$
%asserting that in $V$ and every forcing extension, the Keisler completeness
%theorem holds).

%(Note that this is not comparable to $M\vD \qb P\text{ normal}\qe$,
%since in $M$ there could e.g. be a largest cardinal, so the usual
%definition of normality using the notion of arbitrary 
%large cardinals would not make much sense.)
%(By the way, there will of course always be non-normal candidates as well.)

\begin{Thm}\label{maintheorem}
Assume that $P$ is transitive nep (with respect to a strongly normal $\ZFCx$)
and Borel $I^+$-preserving in $V$
and every forcing extension of $V$.
Then $P$ preserves generics (for unboundedly many candidates)
and therefore $P$ is $I^+$-preserving.
\end{Thm}

We will start with showing that 
for all candidates $M$ and $p\in P\rzM$, $\aGen(M,p)$ 
is nonempty:

\begin{Lem}\label{hilfslemma1}
If $P$ is Borel $I^+$-preserving,  $A\in \IpB$, $M$
a candidate and $p\in P\rzM$, then
$\aGen(M,p)\cap A^V\neq\emptyset$.
\end{Lem}

\begin{proof}
Let $G$ be $P$-generic over $M$ and $V$ and contain $p$. 
In $V[G]$, $\Gen(M\vG)$ is of measure 1, and $A^V$ is positive
(since $P$ is Borel $I^+$-preserving). 
So there is an $\es\in \Gen(M\vG)\cap A^V$.
Let $p'\leq p$ force all this (in particular ``$G$ is $P$-generic over $M$'',
so $p'$ is $M$-generic). Then $p'$ witnesses that
$\es\in \aGen(M,p)$.
\end{proof}

Before we proceed, we take a look once more at strongly normal theories, to
make sure that the models we will be using in the proof really are
$\ZFCx$-candidates. Intuitively, the reader can think of ZFC models instead of
$\ZFCx$ (formally that  would require a few inaccessibles) and elementary
submodels of the universe instead of $H(\chi)$ (that would be more complicated
to justify formally).

$\ZFCx$ is strongly normal, so for any forcing notion $R$, $\chi'$ regular and
large, $1_R\forc H(\chi')^{V[G]}\vD\ZFCx$.  For $p\in R\subseteq H(\chi)$,
$\chi'\gg\chi$ regular, $\n\tau\in H(\chi')$, the following are equivalent:
$H(\chi')\vD\qb p\forc_R \varphi(\n\tau)\qe$ and $p\forc_R(H(\chi')^{V[G]}\vD
\varphi(\n\tau))$.  So in $H(\chi')$ the following holds: For all small
forcings $R$, $1_R\forc_R\ZFCx $. 
%In other words, we can (recursively) enlarge $\ZFCx$ to a $\ZFCxx$
%such that the following  holds:\\
%(1) $\ZFCxx$ is still strongly normal\\
%(2) if $M$ is a $\ZFCxx$-model, and $M[G]$ a generic extension of 
%$M$, then $M[G]$ is a $\ZFCx$-model.

``$P$ is  Borel $I^+$-preserving'' is absolute between
$V$ and $H(\chi)$ for $\chi>2^\al0$ regular,
since for every $A\in\IpB\subset H(\chi)$, 
$p \forc_P A^V\in I$ iff $p\forc_P H(\chi)^{V[G_P]}\vD A^V\in I$
iff $H(\chi)\vD p\forc_P A^V\in I$.
Also, ``$P$ is transitive nep'' is absolute:
every countable transitive candidate $M$ and every $p\in P$
is in $H(\chi)$, and $p\forc_P(G_P\cap P^M$ is $M$-generic$)$
is absolute by the same argument.
In the same way we see the following:
If $R\in H(\chi)$, $\chi\ll\chi'$, then
``$\forc_R P\text{ is transitive nep and Borel }I^+\text{-preserving}$''
is absolute between $V$ and $H(\chi')$, and therefore true in 
$H(\chi')$ according to our assumption.

So every forcing extension $M'$ (by a small forcing) of  $H(\chi')$
(or a transitive collapse of an elementary submodel of $H(\chi')$)
as well as $H(\chi)^{M'}$ (for $\chi$ large with respect to the forcing)
will satisfy $\ZFCx$ and think that $P\text{ is
transitive nep and Borel }I^+\text{-preserving}$.

%Since $\ZFCx$ is normal, the theory $\ZFCxx$,
%$\ZFCx$ plus 

%Since $\ZFCx$ is strongly normal, in any universe $V'$, $H(\chi)^{V'}$
%satisfies $\ZFCx$ (if $\chi$ is large enough).
Now we can proceed with the proof of the theorem:
Fix $\chi_1\ll\chi_2\ll\chi_3$ regular such that $H(\chi_i)\vD\ZFCx$.
Let $N\esm H(\chi_3)$ contain $P,\chi_1,\chi_2$. 
Clearly there are unboundedly many
such $N$. Let $M$ be the transitive collapse of $N$.  We want to show that $P$
preserves generics for $M$.

%Assume $P$ is transitive nep with respect to some strongly normal $\ZFCx$.  Let
%$T$  be the theory $\ZFCx$ together with the sentence ``$\ex \chi_1 \ll
%\chi_2\text{ regular such that }H(\chi_i)\vD\ZFCx$''.  Clearly $T$ is strongly
%normal as well.  (If $H(\chi)\vD\ZFCx$ for all regular cardinals $\chi\geq
%\chi_1$, and $\chi_1\ll \chi_2\ll\chi_3$ are all regular, then for all
%$\chi\geq \chi_3$ regular, $H(\chi)\vD T$.) So there is a strongly normal
%theory $\ZFCxx\supseteq T$ such that for every $\ZFCxx$-candidate $M$, and
%every generic extension $M[G]$ of $M$, $M[G]\vD T$.  Since $\ZFCxx\supseteq
%\ZFCx$, $P$ is transitive nep with respect to $\ZFCxx$ as well.
%Assume $M$ is a $\ZFCxx$-candidate. 

In $M$,
%$2^\al0<\chi_1$, $2^{\chi_1}<\chi_2$, 
let $H_1\DEFEQ H(\chi_1)\vD\ZFCx$.
Let $R_i$ (in $M$) be the collapse of
$H(\chi_i)$ to $\omega$. (I.e. $R_i$ consists of finite functions
from $\omega$ to $H(\chi_i)$.)
Let $\es\in\Gen(M)$, $p_0\in P\rzM$.
We have to show that $\es\in\aGen(M,p_0)$.
Let $G_Q\in V$ be an $M$-generic filter such that
$\h[G_Q]=\es$, and let
$G_R\in V$ be $R_2$-generic over $M[G_Q]$, $M'=M[G_Q][G_R]$.

\begin{Lem}\label{hilfslemma2} 
$M'\vD\qb H_1 \text{ is a $\ZFCx$-candidate, }\es\in\aGen(H_1,p_0)\qe$.
\end{Lem}

If this is correct, then theorem \ref{maintheorem} follows:
Assume 
$M'\vD\qb 
p'\leq  p_0\ H_1\text{-generic},\ 
p'\forc\es\in\Gen(H_1[G_P])\qe$.
$M'$ is a 
$\ZFCx$-candidate, 
so we can find a $p''\leq p'$ be $M'$-generic. Then $p''$ is $H_1$ generic
and therefore $M$ generic as well
(since $\pow(P)\cap M=\pow(P)\cap H_1$), 
and $p''\forc \es\in\Gen(M[G_P])$.

\begin{proof}[Proof of lemma \ref{hilfslemma2}] 
It is clear that 
$H_1$ is a $\ZFCx$-candidate in $M'$.
Assume towards a contradiction, that 
$M'\vD\qb \es\notin \aGen(H_1,p_0)\qe$. Then this is 
forced by some $q\in G_Q$
and $r\in R_2$, but since $R_2$ is homogeneous, 
without loss of generality
 $r=1$, i.e.\begin{equation}\tag{$\ast$}
M\vD\qb q\forc_Q\, \forc_{R_2} \es\notin \aGen(H_1,p_0)\qe.\end{equation}

\begin{table}[tb]
\newcommand{\areins}{\ruDashto^{R_1/Q}_{\tilde{G}_1}}
\newcommand{\arzwei}{\rDashto^{R_2}_{\tilde{G}_1\ast\tilde{G}_2}}
\begin{diagram}
M &           & \rTo^{R_1}&       &M[G_{R_1}]\EQDEF M_1 & %\vD\es\in\aGen(H_1)
      & \\
  &\rdDashto^Q&           &\areins&                     & \rdTo^{R'}_{\tilde{G}_2} & \\
  &           &M[\eot]    &       &\arzwei              &                          & M[\eot][G_{R_2}]=M_2  \\
\end{diagram}
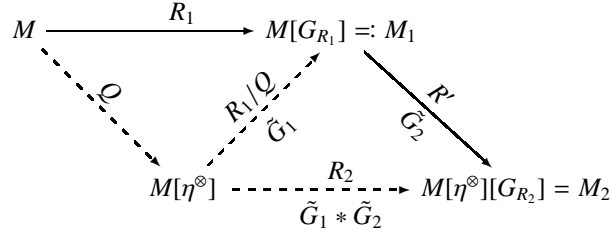
\caption{The models used in the proof of \ref{hilfslemma2}}
\label{table:1}
\end{table}
Now we are going to construct the models of table \ref{table:1}:
First, choose a $G_{R_1}\in V$ which is 
$R_1$-generic over $M$, and let $M_1=M[G_{R_1}]$.
In $M_1$, pick $\eot\in\aGen(H_1,p_0)\cap B^M_q$. 
(We can do that by lemma \ref{hilfslemma1}, since we know that 
$P$ is Borel $I^+$-preserving in $M_1$). 
Since $\aGen\subseteq \Gen$,
$M_1\vD\qb \ex\, G_Q^\otimes\ Q\text{-generic over }H_1$ such that $q\in G_Q^\otimes,
\h[G_Q^\otimes]=\eot \qe$.
This $G_Q^\otimes$ clearly is $M$-generic as well (since 
$M\cap\pow(Q)=H_1\cap\pow(Q)$),
so we can factorize $R_1$ as $R_1=Q\ast {R_1/Q}$ such that
$G_{R_1}=G_Q^\otimes\ast \tilde{G}_1$. 

Now we look at the forcing $R_2=R_2^M$ in $M[\eot]=M[G_Q^\otimes]$.  
$R_2$ forces that
$R_1$ is countable and therefore equivalent to Cohen forcing.
$R_1/Q$ is a subforcing of $R_1$. Also,
$R_2$ adds a Cohen real. So $R_2$ can be factorized as $R_2=(R_1/Q)\ast R'$,
where $R'=(R_2/(R_1/Q))$. We already have $\tilde{G}_1$,
a $(R_1/Q)$-generic filter
over $M[G_Q^\otimes]$, now choose $\tilde{G}_2\in V$ $R'$-generic over $M_1$,
and let $G_{R_2}=\tilde{G}_1\ast\tilde{G}_2$ So $G_{R_2}\in V$ is
$R_2$-generic over $M[G_Q^\otimes]$, $M_2\DEFEQ M[\eot][G_{R_2}]$.

Let $H_2$ be $H(\chi_2)^{M_1}$. 
$H_2\vD\ZFCx$.
Also,
$H_2\vD\qb p_1\leq p_0\text{ is }H_1\text{-generic}, p_1\forc \eot\in\Gen(H_1[G_P])\qe$ (since this is absolute between the universe $M_1$ and
$H_2=H(\chi_2)^{M_1}$).
In $M_2$, $H_2$ is a 
$\ZFCx$-candidate. Let in $M_2$, $p_2\leq p_1$ be $H_2$-generic.
Then (in $M_2$), $p_2$ witnesses that $\es\in\aGen(H_1,p_0)$, 
a contradiction to ($\ast$).
\end{proof}

%\section{Bibliography}
\bibliographystyle{plain}
\bibliography{logik,listb,listx}

%\printindex
\end{document}